\documentclass[graybox]{svmult}

\usepackage{mathptmx}       
\usepackage{helvet}         
\usepackage{courier}        
\usepackage{type1cm}        
\usepackage{makeidx}         
\usepackage{graphicx}        
\usepackage{multicol}        
\usepackage[bottom]{footmisc}

\usepackage{amsmath,amstext,amssymb,amsxtra}

\usepackage{mathtools}
\mathtoolsset{showonlyrefs,showmanualtags} 

\usepackage{enumerate}
\usepackage{pgf,tikz}
\usetikzlibrary{arrows,positioning}

\usepackage{multirow}

\makeindex             

\begin{document}

\title*{The Supremum Norm of the Discrepancy Function: Recent Results and Connections}
\titlerunning{Sup-Norm of the Discrepancy Function}
\author{Dmitriy Bilyk and Michael Lacey}
\institute{Dmitriy Bilyk \at School of Mathematics, University of Minnesota, Minneapolis MN 55455, USA, \email{dbilyk@math.umn.edu}
\and Michael Lacey \at School of Mathematics, Georgia Institute of Technology, Atlanta GA 30332, USA, \email{lacey@math.gatech.edu}}

%

%
\maketitle

\abstract{A great challenge in the analysis of the discrepancy function  $ D_N$ is to obtain universal lower bounds 
on the $ L ^{\infty }$ norm of $ D_N$ in dimensions $ d \ge 3$.  
It follows from the $ L ^{2}$ bound of Klaus Roth that $ \lVert D_N\rVert_{\infty } \ge \lVert D_N\rVert_{2} 
\gtrsim (\log N) ^{(d-1)/2}$.  It is conjectured that the $ L ^{\infty }$ bound is significantly larger, but the only 
definitive result is that of Wolfgang Schmidt in dimension  $d=2$.    Partial  improvements of the 
Roth exponent $ (d-1)/2$ in higher dimensions have been  established  by the authors  and Armen Vagharshakyan.  We survey these results, 
 the underlying methods, and some of their connections to other subjects in probability, approximation theory, and analysis.
}

\section{Introduction} 

We survey recent results on the sup-norm of the discrepancy function. 
For integers $ d\ge 2$, and $ N\ge 1$, let $ \mathcal P_N  \subset  [0,1] ^{d} $ be a finite point set  with cardinality $  \sharp  \mathcal P_N=N $. Define the associated discrepancy function by 
\begin{equation}  \label{e.discrep}
D_N(x)=\sharp( \mathcal P_N \cap [ 0, x))-N \lvert  [ 0, x) \rvert,
\end{equation}
where  $ x= (x_1 ,\dotsc, x_d)$ and $[ 0, x)=\prod _{j=1}^d [0,x_j)$ is a rectangle with antipodal corners at  $ 0$ and $ x$, and $|\cdot |$ stands for the $d$-dimensional Lebesgue measure.
The dependence upon the selection of points $\mathcal P_N$ is suppressed, as we are mostly interested  in bounds 
that are universal in $ \mathcal P_N$.  The discrepancy function $D_N$ measures equidistribution of $\mathcal P_N$: a  set of points is \emph{well-distributed} if  this function is small in some appropriate function space.

It is a basic fact of the theory of irregularities of distribution that relevant norms of this function  in dimensions $2$ and higher must tend to infinity as $N$ grows. 
The classic results are due to   Roth \cite{MR0066435} in the case of the $ L ^{2}$ norm
and  Schmidt \cite{MR0491574} for   $L ^{p} $, $ 1< p < 2$.  

\begin{theorem}\label{t.DP} 
For $1<p<\infty$ and any collection of points $ \mathcal P_N\subset [0,1] ^{d}$, we have
\begin{equation}\label{e.DP}
\lVert  D_N \rVert _{ p} \gtrsim (\log N) ^{(d-1)/2}\,. 
\end{equation}
Moreover, we have the endpoint estimate
\begin{equation}\label{e.DPend}
\lVert  D_N \rVert _{ L (\log L )  ^{(d-1)/2}} 
\gtrsim (\log N) ^{(d-1)/2}.
\end{equation}
\end{theorem}

The symbol ``$\gtrsim$" in this paper stands for ``greater than a constant multiple of", and the implied constant may depend on the dimension, the function space, but {\it{not}} on the configuration $\mathcal P_N$ or the number of points $N$. The Orlicz space notation, such as $L (\log L)^\beta$, is explained in the next section, see \eqref{orlicz}. 

We should mention that there exist sets $ \mathcal P_N$ that meet the $ L ^{p}$  bounds \eqref{e.DP} in all dimensions.  
This remarkable fact is established by beautiful and quite non-trivial constructions of the point distributions $ \mathcal P_N$. 
We refer to the reader to one of the very good references \cite{MR903025,MR2683394,MR1470456} on the subject for 
more information on this important complement to the subject of this note.

While the previous theorem is quite adequate for $ L ^{p}$, $1<p<\infty$, the endpoint cases 
of $ L ^{\infty }$ and $ L ^{1}$ are not amenable to the same techniques. Indeed, the extremal $ L ^{\infty }$ bound should be larger than the average $L^2$ norm. In dimension $d=2$ the endpoint estimates are known --    it is the theorem of Schmidt \cite{MR0319933}. 

\begin{theorem}  \label{t.schmidt}
The following estimate is valid for all collections $\mathcal P_N\subset [0,1]^2 $:  
\begin{align}  \label{e.schmidt}
\lVert  D_N  \rVert _{\infty} {}\gtrsim{} \log N . 
\end{align}
\end{theorem} 

This is  larger than Roth's $L^2$ bound by $ \sqrt {\log N}$.  
The difference between the two estimates points to the fact that for extremal choices of sets $ \mathcal P_N$,  the $ L ^{\infty }$ norm of $ D_N$ is obtained on a set so small it cannot be seen on the scale of $L ^{p} $ spaces.  We will return to this point below.

In dimensions $3$ and higher partial results began with a breakthrough work 
of J.~Beck \cite{MR1032337} in dimension $ d=3$.  
The following result is  due to Bilyk and Lacey \cite{MR2414745} in dimension $ d=3$, 
and Bilyk, Lacey, Vagharshakyan \cite{MR2409170} in dimensions $ d\ge 4$.  

\begin{theorem}\label{t.beck} 
In dimensions $ d\ge 3$ there exists $\eta = \eta (d)\ge c / d ^2 $  for which the following estimate 
holds for all collections $\mathcal P_N\subset [0,1]^d $:  
\begin{align}  \label{e.beck}
\lVert  D_N  \rVert _{\infty} {}\gtrsim{} (\log N) ^{ \frac{d-1}{2}+ \eta }\,. 
\end{align}
\end{theorem}

This is  larger than Roth's bound by $ (\log N) ^{\eta }$.  
Beck's original result   in dimension $ d=3$ had a  much smaller 
doubly logarithmic term $(\log \log N)^{\frac18 -\varepsilon}$ in place of   $ (\log N) ^{\eta }$.  The proof strategy begins with the fundamental orthogonal function method 
of Roth and Schmidt, which we turn to in the next section.  In \S\ref{s.small} we turn to  a closely related 
combinatorial inequality for ``hyperbolic" sums of  multiparameter Haar functions. It serves as the core question which has related the progress on 
lower bounds for the discrepancy function to questions in probability and approximation theory.  
Based upon this inequality, it is natural to conjecture that the optimal form of  the $ L ^{\infty }$ estimate is 

\begin{conjecture}\label{j.D} 
In dimensions $ d\ge 3$  there holds $ \lVert D_N\rVert_{\infty } \gtrsim (\log N) ^{d/2}$. 
\end{conjecture}

We should mention that at the  present time  there is no consensus among the experts about the sharp form of the conjecture (in fact, a great number of specialist believes that $\lVert D_N\rVert_{\infty } \gtrsim (\log N) ^{d-1}$ is the optimal bound, which is supported by the best known examples). However,  in this paper we shall advocate our belief in Conjecture \ref{j.D} by comparing it to other sharp conjectures in various fields of mathematics. In particular, the sharpness of Conjecture \ref{j.small} in \S\ref{s.small} suggests that  the estimate above is the best that could be obtained by the orthogonal function techniques.

The reader can consult the papers \cite{MR2414745,MR2409170}, as well as the surveys of 
the first author \cite{MR2817765,survey} for more detailed information.  {{This research is supported in part by NSF grants  DMS 1101519, 1260516 (Dmitriy Bilyk),  DMS 0968499, and  a grant from the Simons Foundation \#229596 (Michael Lacey).}}

\section{The Orthogonal Function Method} 

All progress on these universal lower bounds has been based upon  the orthogonal function method, initiated by Roth, with the modifications of Schmidt, as 
presented here.  
Denote the family of all dyadic intervals $ I\subset [0,1]$ by $ \mathcal D$.  Each dyadic interval 
$ I$ is the union of two dyadic intervals $ I _{-}$ and $I_+$, each  of exactly half the length of $ I$, representing
the left and right halves of $ I$ respectively.  Define the Haar function associated to $ I$ 
by $ h_I = - \chi _{I _{-}}+ \chi _{I _{+}}$.  Here  and throughout we will use the $ L ^{\infty }$ (rather than $L^2$)
normalization of the  Haar functions.  

In dimension $ d$, the $ d$-fold product $ \mathcal D ^{d}$ is the collection of dyadic intervals in $ [0,1] ^{d}$. 
Given $ R= R_1 \times \cdots \times R_d \in \mathcal D ^{d}$, the Haar function associated with $ R$ is the 
tensor product 
\begin{equation*}
h _{R} (x_1 ,\dotsc, x_d) = \prod _{j=1} ^{d} h _{R_j} (x_j) \,. 
\end{equation*}
These functions are pairwise orthogonal as $ R\in \mathcal D ^{d}$ varies.  

For a $d$-dimensional vector  $ r= (r_1 ,\dotsc, r_d)$ with non-negative integer coordinates 
let $ \mathcal D _{r}$ be the set of those $ R\in \mathcal D ^{d}$ that 
for each coordinate $ 1\le j \le d$, we have $ \lvert  R _{j}\rvert = 2 ^{- r_j} $. These rectangles partition $ [0,1] ^{d}$.  
We call $ f _{r}$ an $ r$-function (a generalized Rademacher function) if for some choice of signs $ \{\varepsilon _{R} \;:\; R\in \mathcal D _{r}\}$, we have 
\begin{equation*}
f _{r} (x) = 
\sum_{R \in \mathcal D _{r}} \varepsilon _R h _{R} (x)\,. 
\end{equation*}

The following is the crucial lemma of the method.  Given an  integer $ N$, we set $ n = \lceil 1 + \log_2 N \rceil$, where $\lceil x\rceil$ denotes the smallest integer greater than or equal to $x$. 

\begin{lemma}\label{l.r}  
In all dimensions $ d\ge 2$ there is a constant $ c_d > 0 $ such that 
for each $r$ with  $\lvert  r\rvert := \sum_{j=1} ^{d} r_j=n$, there is an $ r$-function $ f _{r}$ with 
$ \langle D_N, f _{r} \rangle \ge c_d$.   
Moreover, for all $ r$-functions there holds  $ \lvert  \langle D_N , f _{r} \rangle\rvert \lesssim N 2 ^{- \lvert  r\rvert }$.  
\end{lemma}

The proof of the lemma is straightforward, see e.g. \cite{MR0066435,MR0491574,MR2817765}.  With this lemma at hand, the proof of Roth's Theorem in $ L ^{2}$ is as follows.  
Note that the requirement that $ \lvert  r\rvert=n $ says that the coordinates of $ r$ must partition $ n$ into 
$ d$ parts. It follows that the number  of  ways to select the coordinates of $ r$ is bounded above and below by a multiple of $  n ^{d-1}$, agreeing with the simple logic that there are $d-1$ ``free" parameters: $d$ dimensions minus the  restriction  $| r |  = n$. Set $ F _{d} =  \sum_{r \;:\; \lvert  r\rvert=n }  f _{r} $. Orthogonality implies that $\|F_d\|_2 \lesssim n^{(d-1)/2}$. Hence, by Cauchy--Schwarz
\begin{align}
 n ^{d-1} & \lesssim \sum_{r \;:\; \lvert  r\rvert=n } \langle D_N , f _{r} \rangle = \langle D_N , F_d \rangle
\\ \label{e.middle}
& \le \lVert D_N\rVert_{2} \cdot \|F_d\|_2
 \le  \lVert D_N\rVert_{2}  \cdot n ^{ (d-1)/2} \,. 
\end{align}
The universal lower bound $ n ^{ (d-1)/2} \lesssim \lVert D_N\rVert_{2}$ follows.

Deeper properties of the discrepancy function  may be deduced from finer properties of $ r$-functions.   A key property is the classical Littlewood--Paley inequality  for Haar functions:

\begin{theorem}\label{t.LP} For $p\ge 2$, we have the inequality 
\begin{equation} \label{e.LP}
\Bigl\lVert  \sum_{I\in \mathcal D} \alpha_I h_I  \Bigr\rVert_{p} 
\le C \sqrt p\, \Bigl\lVert  \Bigl[
\sum_{I\in \mathcal D} \frac { \lvert  \alpha_I\rvert ^2 } {\lvert  I\rvert ^2  }  \chi _I 
\Bigr] ^{1/2}  \Bigr\rVert_{p} \,,
\end{equation}
where $ C$ is an absolute constant, and the coefficients $ \alpha _I$ take values in a Hilbert space $ \mathbf H$.  
\end{theorem}

The right-hand side is the Littlewood--Paley (martingale) square function of the left hand side.  This inequality can be viewed as an extension of orthogonality and Parseval's identity to values of $p$ other than $2$, and it is often useful to keep track of the growth of $ L ^{p}$ norms.  The fact that one can allow Hilbert space value 
coefficients permits repeated application of the inequality.  The role of the Hilbert space valued coefficients is the focus of \cite{MR1439553}, 
which includes more information about  multiparameter harmonic analysis, relevant to this subject.

Consider the dual function in \eqref{e.middle}, $ F _{d} =  \sum_{r \;:\; \lvert  r\rvert=n }  f _{r} $.
As discussed earlier, the index set  $ \{r \;:\; \lvert  r\rvert=n \}$ has $ d-1$ free parameters.  
The function $ F_d$ is a Haar series in the first variable, so the inequality  \eqref{e.LP} applies.  
On the right-hand side, the square function can be viewed as an $ \ell ^2 $-valued Haar series in the second variable, hence \eqref{e.LP}  applies again, see \cite{MR2409170,MR2817765} for details.   Continuing this $ d-1$ times, one arrives at 
\begin{equation} \label{e.slp}
\lVert F_d\rVert_{p} \lesssim p ^{ (d-1)/2} n ^{ (d-1)/2} \,, \qquad 2\le  p < \infty \,. 
\end{equation}
Repeating \eqref{e.middle} verbatim (with H\"{o}lder replacing Cauchy--Schwarz), one obtains $ n ^{ (d-1)/2} \lesssim \lVert D_N\rVert_{q} $ for $ 1< q < 2$.  

If one is interested in endpoint estimates, it is useful to rephrase the  inequalities for $ F_d$ above in the language of Orlicz spaces.  For a convex increasing function $\psi: \mathbb R_+ \rightarrow \mathbb R_+$ with $\psi(0) = 0$, the Orlicz space $L^\psi$ is defined as the space of measurable functions $f: [0,1]^d \rightarrow \mathbb R$ for which \begin{equation}\label{orlicz} 
\| f \|_{L^\psi} = \inf \bigg\{ K>0: \int_{[0,1]^d} \psi \big( |f(x)|/K \big) \, dx \le 1 \bigg\}.
\end{equation}
In particular, for $\psi (t) = t^p$ one obtains the standard $L^p$ spaces, while $\operatorname{exp} (L^\alpha)$ and $L (\log L)^\beta$ denote Orlicz spaces generated by functions equal to $e^{t^\alpha}$ and $t \log^\beta t$ respectively, when $t$ is large enough. These spaces serve as refinements of the endpoints of the $L^p$ scale, as for each $1<p<\infty$, $\alpha$, $\beta>0$  we have the embeddings $L^\infty \subset \operatorname{exp} (L^\alpha) \subset L^p$ and $L^p \subset L (\log L)^\beta \subset L^1$.

The polynomial growth in the $ L ^{p}$ norms of $F_d$ \eqref{e.slp} translates into exponential integrability estimates, 
 namely $ \lVert F_d \rVert_{\operatorname {exp}(L ^{2/ (d-1)})} 
\lesssim n ^{ (d-1)/2}$, since 
\begin{equation*}
\| f \|_{\operatorname{exp}(L ^{\alpha })} \simeq
\sup _{p>1} p ^{-1/\alpha } \| f \|_p \,, \qquad \alpha >0 \,. 
\end{equation*}
The dual space to $ {\operatorname {exp}(L ^{2/ (d-1)})}$ is $ L (\log L) ^{ (d-1)/2}$, 
hence we see that 
\begin{equation*}
n ^{ (d-1)/2} \lesssim \lVert D_N\rVert_{L(\log L) ^{(d-1)/2} }  \,. 
\end{equation*}

A well-known result of  Hal{\'a}sz \cite{MR637361} is a  `$ \sqrt {\log L}$' improvement of this estimate in dimension $ d=2$.  
Indeed, we have the following theorem  valid for all dimensions, see \cite{MR2419612}.  

\begin{theorem}\label{t.l} For dimensions $ d\ge 2$, there holds $ \lVert D_N\rVert_{L(\log L) ^{(d-2)/2}} \gtrsim (\log N)^{ (d-1)/2}$. 
\end{theorem}
\noindent Notice that for $d=2$ one recovers Hal\'{a}sz's $L^1$ bound 
\begin{equation}\label{e.L1}
 \sqrt{\log N} \lesssim \| D_N \|_1 .
\end{equation}
  In dimension $ d=2$, the argument of Hal{\'a}sz  can be rephrased into the estimate 
\begin{equation} \label{e.halaszDual}
\sqrt n \lesssim \langle  D_N ,  \sin \bigl(\tfrac c {\sqrt n} F_2 \bigr)\rangle\,, \qquad 0 < c < 1 \textup{ sufficiently small.} 
\end{equation}
This immediately shows that $ \lVert D_N\rVert_{1} \gtrsim \sqrt n$ in dimension $ d=2$.  
There is a relevant endpoint estimate of the Littlewood--Paley inequalities, namely the Chang--Wilson--Wolff inequality 
\cite{MR800004}.  Employing extensions of this inequality and the estimate above, one can give a proof of Theorem \ref{t.l} 
in dimensions $ d\ge 3$.

It is a well-known conjecture that in all dimensions $ d \ge 3$ one has  the estimate \begin{equation}  \lVert D_N\rVert_{1} \gtrsim (\log N)^{ (d-1)/2} \end{equation} on 
the $ L ^{1}$ norm of the discrepancy function. Any improvement of Theorem \ref{t.l} \ would yield progress on this conjecture.

\section{The Small Ball Inequality} \label{s.small}

Lower bounds on the discrepancy function are related through proof techniques to subjects in different areas of mathematics. 
They include, in particular, the    so-called \emph{small deviation inequalities} for the Brownian sheet  in probability theory, 
complexity bounds for certain Sobolev spaces  in approximation theory, and
a combinatorial inequality involving multivariate Haar functions 
in the unit cube.  We refer the reader to the references \cite{MR2414745,MR2409170,MR2817765,survey} for more  information, and emphasize 
that the questions in probability and approximation theory are parts of very broad areas of investigation with additional 
points of contact with discrepancy theory  and many variations of the underlying themes.  

\smallskip 
According to the idea introduced in the previous section,  the behavior of $D_N$ is essentially defined by its projection onto the span of $\{ h_R:\, R\in \mathcal D^d, \,|R|=2^{-n} \}$. It is therefore reasonable to model the discrepancy bounds by estimates of the linear combinations of such Haar functions (we call such sums ``hyperbolic").  The problem of obtaining  lower bounds  for sums of Haar functions supported by rectangles of fixed volume -- known as the {\it Small Ball inequality} -- arises naturally in the aforementioned problems in probability and approximation theory.  While in the latter fields versions of this inequality have important formal implications, its connection to discrepancy estimates is still only intuitive and is not fully understood. However, most known proof methods are easily transferred from one problem to another. The conjectured form of the inequality is the following.
  
\begin{conjecture}\label{j.small}
 \label{small}[The Small Ball Conjecture] For dimension $ d\ge 3$ we have the inequality 
 \begin{equation}\label{e.small}
2 ^{-n} \sum _{\lvert  R \rvert= 2 ^{-n} } \lvert  \alpha_R  \rvert 
{}\lesssim{}  n ^{(d-2)/2  }
\bigg\|  \sum _{R\in \mathcal D^d :\,  \lvert  R \rvert =  2 ^{-n} } \alpha_R h_R \bigg\|_{\infty }
\end{equation}
valid for all real-valued coefficients $ \alpha _R$. 
\end{conjecture}

The subject of the  conjecture is the exact exponent of $ n$ the right-hand side. This conjecture is better, by one square root of $ n$,  than a trivial estimate available from the Cauchy--Schwartz inequality.  Indeed, with $n^{(d-2)/2}$ replaced by $n^{(d-1)/2}$ it holds for the $L^2$ norm:
\begin{align}
\bigg\| \sum_{R\in \mathcal D^d :\, |R|=2^{-n}} \alpha_R h_R\bigg\|_2 & = \bigg( \sum_{|R|=2^{-n}} |\alpha_R|^2 2^{-{n}} \bigg)^\frac12 \label{sbl2}\\
\nonumber & \gtrsim \frac{\sum_{|R|=2^{-n}} |\alpha_R| 2^{-{n}/2}}{ \big( n^{d-1} 2^n \big)^\frac12 } = n^{-\frac{d-1}2} \cdot 2^{-n}  \sum_{|R|=2^{-n}} |\alpha_R|,
\end{align}
where we have used the fact that the total number of rectangles $R\in \mathcal D^d$ is $\approx n^{d-1} 2^n$. This computation is similar in spirit to \eqref{e.middle} establishing Roth's $L^2$ discrepancy bound.
 Generally, the Small Ball Conjecture bears a strong resemblance to Conjecture \ref{j.D} about the discrepancy function. Indeed, in both cases one gains a square root of the logarithm over the $L^2$ bound.

One can  consider a restricted version of inequality \eqref{e.small}, which  appears to contain virtually 
all the complexity of the general inequality and is sufficient for applications:
 \begin{equation}\label{e.signed}
   \Bigl\lVert \sum _{\lvert  R \rvert= 2 ^{-n} } \varepsilon_R h_R  \Bigr\rVert_{\infty } \gtrsim n ^{d/2}\,, \qquad \varepsilon _{R} \in \{-1,0,1\}\,, 
\end{equation}
subject to the requirement that $ \sum _{\lvert  R \rvert= 2 ^{-n} } \lvert  \varepsilon _R\rvert \ge c 2 ^{n} n ^{d-1} $ for a fixed  small constant $ c>0$, in  other words,  at least a fixed proportion of the coefficients $ \varepsilon_R$ are non-zero.
The relation to the discrepancy estimates becomes even more apparent for this form of the inequality. For instance,  the trivial bound \eqref{sbl2} 
becomes 
 \begin{equation}
   \Bigl\lVert \sum _{\lvert  R \rvert= 2 ^{-n} } \varepsilon_R h_R  \Bigr\rVert_{\infty }
   \ge 
    \Bigl\lVert \sum _{\lvert  R \rvert= 2 ^{-n} } \varepsilon_R h_R  \Bigr\rVert_{2 }
   \gtrsim 
   n ^{(d-1)/2}\,.
\end{equation}
Compare this to   Roth's bound \eqref{e.DP}, and compare \eqref{e.signed}  to Conjecture~\ref{j.D}. 
The similarities between the discrepancy estimates and the Small Ball inequality are summarized in Table \ref{table1}. 

A more restrictive version of inequality \eqref{e.small} with $\varepsilon_R =\pm 1$ (the {\it{signed small ball inequality}}) does allow for some proof simplifications, but has no direct consequences.  The papers \cite{MR2375609,MR2409170} study this restricted inequality, 
using only the fundamental inequality -- Lemma~\ref{l.SimpleCoincie} of \S \ref{s.rr}.  
This case will likely continue to be a proving ground for new techniques in this problem.

\begin{table}
\begin{tabular}{| c | c |}
\hline
 &
\\
  {\bf{Discrepancy estimates}} & {\bf{Small Ball inequality (signed)}}
\\[4pt]
\hline
\multicolumn{2}{|c|}{Dimension $d=2$} \\[2pt]
\hline
&
\\[-3pt]
${\|D_N\|_\infty \gtrsim \log N}$  &  $\displaystyle{\bigg\|  \sum_{|R|=2^{-n}} \varepsilon_R h_R \bigg\|_\infty \gtrsim n }$ \,\, 
\\[9pt]
 (Schmidt, '72; Hal\'asz, '81) & (Talagrand, '94; Temlyakov, '95)
\\[4pt]
\hline
\multicolumn{2}{|c|}{Higher dimensions, $L^2$ bounds} \\[2pt]
\hline
 &
\\[-3pt]
{${\| D_N \|_2 \gtrsim ({\log N})^{{(d-1)}/{2}}}$}  &   $\displaystyle{\bigg\|  \sum_{|R|=2^{-n}} \varepsilon_R h_R \bigg\|_2 \gtrsim n^{{(d-1)}/{2}}}$
\\[4pt]
\hline
\multicolumn{2}{|c|}{Higher dimensions, conjecture} \\[2pt]
\hline
 &
\\[-3pt]
{${\| D_N \|_\infty \gtrsim ({\log N})^{{d}/{2}}}$}  &   $\displaystyle{\bigg\|  \sum_{|R|=2^{-n}} \varepsilon_R h_R \bigg\|_\infty \gtrsim n^{{d}/{2}}}$
\\[4pt]
\hline
\multicolumn{2}{|c|}{Higher dimensions, known results} \\[2pt]
\hline
 &
\\[-3pt]
{${\| D_N \|_\infty \gtrsim ({\log N})^{\frac{d-1}{2}+\eta}}$}  &   $\displaystyle{\bigg\|  \sum_{|R|=2^{-n}} \varepsilon_R h_R \bigg\|_\infty \gtrsim n^{\frac{d-1}{2}+\eta}}$
\\[6pt]
 \hline
\end{tabular}\vskip1mm
\caption{Discrepancy estimates and the signed Small Ball inequality} \label{table1}
\end{table}

 Conjecture \ref{j.small} is sharp: for independent random 
selection of coefficients (either random signs or Gaussians), the supremum is  at most $ C n ^{d/2}$,
\begin{equation*}
\mathbb E  \Bigl\lVert \sum _{\lvert  R \rvert= 2 ^{-n} } \alpha_R h_R  \Bigr\rVert_{\infty } \simeq n ^{d/2} \,. 
\end{equation*}
Unfortunately, random selection of coefficients does not seem to be a guide to the sums that are hardest to analyze.  The sharpness of the Small Ball Conjecture justifies our belief in the optimality of  Conjecture \ref{j.D} in discrepancy theory.

\section{Connections to Probability and Approximation Theory}

We briefly touch upon the connections of the Small Ball inequality \eqref{e.small} to problems in other fields. A very detailed account of these relations is contained in \cite{survey}.  \smallskip

 \subsection{Approximation theory:  Metric entropy of classes with dominating mixed smoothness.}

Let $MW ^{p}$  be the image of the unit ball $L^p([0,1]^d)$ under the integration operator $ \big(  \mathcal T f \big)(x) = \int_0^{x_1} ... \int_0^{x_d} f(y)dy$, i.e. in some sense $MW ^{p}$ is the  set  of  functions on $[0,1]^d$ whose mixed derivative
 $\frac{\partial^d \,f}{\partial x_1 \partial x_2 \dots \partial x_d}$ has
 $L^p$ norm bounded by one. This set  is  compact in  the $L^\infty$ metric and its compactness may be
measured by  \emph{covering  numbers}. Let  $N(\varepsilon, p, d) $ be the cardinality of the smallest $\varepsilon$-net of  $MW ^{p} $  in the $L^\infty$ norm. The exact asymptotics of these numbers as $ \varepsilon \downarrow 0$ is a subject of  conjecture.

\begin{conjecture}\label{c.n} For $ d\ge 2 $, we have \,
$\log N(\varepsilon, 2 , d) \simeq  \varepsilon ^{-1}
 (\log 1/\varepsilon)^{d-1/2} \,,$\,\,   as $\varepsilon \downarrow 0$.
\end{conjecture}
\noindent The case $ d=2$  is settled 
\cite{MR95k:60049},
 and the upper bound is known in all dimensions \cite{2000b:60195}. Inequalities similar to the Small Ball Conjecture \eqref{e.small} lead to  lower bounds on the covering numbers. \vskip0.2cm

\subsection{Probability: The small ball problem for the Brownian sheet.}  

Consider the  Brownian sheet $B_d$, i.e.  a centered multiparameter Gaussian process characterized by the
covariance relation $\mathbb E B_d(s) \cdot B_d (t) = \prod_{j=1}^d \min
\{s_j,t_j\}$. The problem deals with the precise behavior of  $\mathbb P (\| B \|_{C([0,1]^d)} < \varepsilon ) $, the {\it{small deviation (or small ball) probabilities}} of $B_d$.

There is  an exciting formal equivalence established by Kuelbs and Li  \cite{MR1184915,MR1237989} 
between the small ball probabilities  and the metric entropy of the unit ball of the
reproducing kernel Hilbert space, 
which in the case of the Brownian sheet is $WM^2$. This yields an equivalent conjecture:


\begin{conjecture}\label{c.sheet}
In dimensions $d\ge 2$, for the Brownian sheet $B$ we have \,\,$$-\log
\mathbb P (\| B\|_{C([0,1]^d)} < \varepsilon ) \simeq
{\varepsilon}^{-2} ( \log 1/\varepsilon )^{2d-1}, \quad \varepsilon
\downarrow 0 .$$
\end{conjecture}
The upper bounds are known for $d\ge 2$
 \cite{2000b:60195}, while the lower bound for $d=2$ has been obtained by Talagrand \cite{MR95k:60049} using \eqref{e.small}.  It is worth mentioning that Conjecture \ref{c.sheet} explains the nomenclature {\it small ball inequality}.

 \subsection{Summary of the connections}
 
The connections between the Small Ball Conjecture and these problems  is illustrated in Figure 1. Solid arrows represent  known formal implications, while a dashed line denotes an informal heuristic relation. Hopefully, other lines, as well as other nodes, will be added to this diagram in the future. In particular, we expect that the theory of empirical processes may connect the discrepancy bounds to the small deviation probabilities. 

 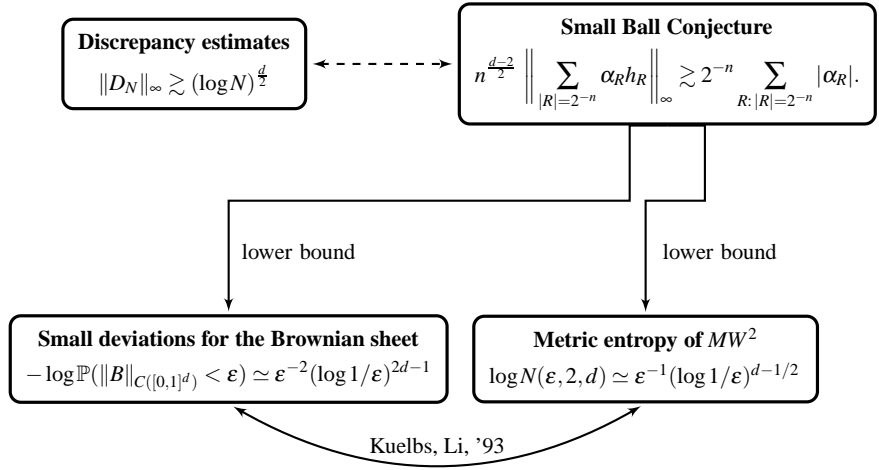
\begin{figure}
\centering
\begin{tikzpicture}[node distance=1cm, auto] 
\tikzset{
    mynode/.style={rectangle,rounded corners,draw=black, top color=white, bottom color=white,very thick, inner sep=0.5em, minimum size=3em, text centered},
    myarrow/.style={->, >=latex', shorten >=1pt, thick},
    mylabel/.style={text width=7em, text centered} 
}  
\node[mynode] (SBI) {\begin{tabular}{c} \textbf{Small Ball Conjecture}\\[1mm]  $\displaystyle{n^{\frac{d-2}{2}}\,\, \bigg\| \sum_{ |R|=2^{-n}} \alpha_R h_R\bigg\|_\infty
\gtrsim 2^{-n}  \sum_{R:\, |R|=2^{-n}} |\alpha_R|.}$ \end{tabular}};
\node[mynode, left=2cm of SBI] (Disc) {\begin{tabular}{c} \textbf{Discrepancy estimates}\\[1mm]  $\displaystyle{\| D_N \|_\infty \gtrsim (\log N)^{\frac{d}2}}$ \end{tabular}};
\node[below=3cm of SBI] (dummy) {}; 
\node[mynode, left=2.8cm of dummy] (BrSh) {\begin{tabular}{c} \textbf{Small deviations for the Brownian sheet} \\[1mm] $-\log \mathbb P (\| B\|_{C([0,1]^d)} < \varepsilon ) \simeq {\varepsilon}^{-2} ( \log 1/\varepsilon )^{2d-1}$ \end{tabular}};  
\node[mynode, right=-2.7cm of dummy] (MetEn) {\begin{tabular}{c}  \textbf{Metric entropy of $MW^2$}\\[1mm]    $\log N(\varepsilon, 2 , d) \simeq  \varepsilon ^{-1}
 (\log 1/\varepsilon)^{d-1/2} \,$ \end{tabular}};
\node at (-4.9,-2.5)  [mylabel] (label1) {lower bound};  
\node at (0.7,-2.5)  [mylabel] (label2) {lower bound};  

\draw[myarrow] (SBI.south) -- ++(-.5,0) --  ++(0,-1) -|    (BrSh.north);	
\draw[myarrow] (SBI.south) -- ++(.5,0) -- ++(0,-1) -|  (MetEn.north);
 
\draw[<->, >=latex', shorten >=2pt, shorten <=2pt, bend right=30, thick] 
    (BrSh.south) to node[auto] {Kuelbs, Li, '93}(MetEn.south); 

\draw[<->, >=latex', shorten >=2pt, shorten <=2pt, dashed, thick] 
    (Disc.east) to node {} (SBI.west); 

\end{tikzpicture} 
\medskip
\caption{Connections between the Small Ball Conjecture and other problems} \label{fig1}
\end{figure}

\section{Riesz Product Techniques}

The only case in which the Small Ball  inequality  \eqref{e.small} is known in full generality is dimension $ d=2$,  which was proved by M. Talagrand \cite{MR95k:60049}. 

\begin{theorem}\label{t.tala} In dimension $ d=2$, there holds for all $ n$, 
\begin{equation*}
2 ^{-n} \sum _{\lvert  R \rvert= 2 ^{-n} } \lvert  \alpha_R  \rvert 
{}\lesssim{}  
\Bigl\lVert  \sum _{\lvert  R \rvert\ge 2 ^{-n} } \alpha_R h_R \Bigr\rVert_{\infty } .
\end{equation*}
\end{theorem}
 
Soon after M. Talagrand  proved Conjecture \ref{j.small} in dimension $d=2$, V. Temlyakov \cite{MR96c:41052}  has given an alternative elegant  proof of this inequality, which  strongly  resonated with  the argument of Hal\'{a}sz \cite{MR637361} for  \eqref{e.schmidt}.  We shall present this technically simpler proof and then explain the adjustments needed to obtain the discrepancy bound.

All the endpoint estimates in dimension $ d=2$ are based upon a  very special property of the two-dimensional Haar functions and  the associated $ r$-functions, {\bf product rule}: if $R,\, R' \in \mathcal D^2$ are not disjoint, $R\neq R'$,  and  $|R| = |R'|$, then 
\begin{equation}\label{productrule}
h_R \cdot h_{R'} = \pm h_{R \cap R'},
\end{equation}
i.e. the product of two Haar functions is again Haar, or equivalently,  if $ \lvert  r\rvert= \lvert  s\rvert =n  $, 
then the product $ f _{r} \cdot f _{s} = f _{t} $ is also an $ r$ function, where $ t = (\min \{ r_1,  s_1 \}, \min\{ r_2 , s_2\} )$.  In higher dimensions two different boxes of the same volume may coincide in one of the coordinates, in which case  $h_{R_k} \cdot h_{R'_k} = h_{R_k}^2 = \mathbf 1_{R_k}$. This loss of orthogonality leads to major complications in dimensions three and above.

\begin{proof}
 For each $j =0,\dots , n$ consider the $r$-functions $\displaystyle{f_{(j,n-j)} = \sum_{\substack{|R|=2^{-n},\,\\ |R_1|=2^{-j}}} \textup{sgn} (\alpha_R) h_R}$.  In    dimension $d=2$ the summation conditions   uniquely define the shape of a dyadic rectangle.   The product rule drives this argument.  We construct the following Riesz product
\begin{equation}\label{psi}
\Psi := \prod_{j=1}^n \bigg( 1 + f_{(j,n-j)} \bigg) =1 + \sum_{R\in \mathcal D^d :\, |R|=2^{-n}} \textup{sgn} (\alpha_R) h_R\, + \,\, \Psi_{> n},
\end{equation}
where, by the product rule, $\Psi_{> n}$ is a linear combination of Haar functions supported by rectangles of area less than $2^{-n}$, and make three simple observations
\begin{enumerate}[(i)]
\item  $\Psi \ge 0 $, since each factor is  either $0$ or $2$. 
\item Next,  $\int \Psi(x) dx = 1$. Indeed, expand the product in \eqref{psi} -- the initial term is  $1$, while all the higher-order terms are    Haar functions with mean zero. 
\item Therefore $\Psi$   has $L^1$ norm $1$: $ \|\Psi\|_1=1$.
\end{enumerate}

\noindent By the same token, using orthogonality,
\begin{align}
\label{sbproof1} \bigg\| \sum_{|R|=2^{-n}} \alpha_R h_R \bigg\|_\infty  &\ge \bigg\langle \sum_{|R|=2^{-n}} \alpha_R h_R, \Psi \bigg\rangle 
  = \,\, 2^{-n} \cdot \sum_{|R|=2^{-n}} |\alpha_R|,
\end{align}
since $\langle h_R, h_R \rangle = 2^{-n}$.
\end{proof}

 Rather than proving Schmidt's discrepancy lower bound, we shall  explain how the above argument could be adapted to obtain Hal\'{a}sz's  proof of \eqref{e.schmidt}. These are the necessary changes:
\begin{itemize}
\item {\it Building blocks:} Instead of the $ r$-functions $f_{(j,n-j)} = \sum  \textup{sgn} (\alpha_R) h_R$ used above, we take the  $ r$-functions   provided by Lemma \ref{l.r} with the property that $\langle D_N, f_r \rangle \gtrsim 1$. 
\item {\it Riesz product:} The test function  $\Psi := \prod_{j=1}^n \bigg( 1 + f_{(j,n-j)} \bigg)$ should be replaced by a slightly more complicated $\Phi =  \prod_{j=1}^n \bigg( 1 + \gamma f_{(j,n-j)} \bigg) - 1$, where $\gamma >0$ is a  small constant.
\end{itemize}
These adjustments play the following roles: $-1$ in the end forces  the ``zero-order" term $\int D_N (x) dx$ to disappear, while a suitable choice of the small constant $\gamma$ takes care of the ``higher-order"  terms and ensures that their contribution is small.  Otherwise, the proof of \eqref{e.schmidt} is verbatim the same as the proof of the two-dimensional Small Ball Conjecture; the details can be found in \cite{MR637361,MR2683232,MR2409170,MR2817765} etc. The Small Ball Conjecture may therefore be viewed as a linear term in the discrepancy estimates.  These same comments apply to the proof of the $ L ^{1}$ estimate \eqref{e.halaszDual} of Hal{\'a}sz.


\begin{table}
\begin{tabular}{| c | c  |}
\hline
&
\\
{\bf{Discrepancy function}} & {\bf{Lacunary Fourier series}}
\\[4pt]
$D_N(x) = \# \{ \mathcal P_N \cap [0,x) \} - Nx_1 x_2$ &   $f(x) \sim \sum_{k=1}^\infty c_k \sin n_k x $, \,\,$\frac{n_{k+1}}{n_k} >\lambda>1$
\\[5pt]
\hline
&
\\[-3pt]       
$\textcolor{red}{\| D_N \|_2 \gtrsim \sqrt{\log N}}$    & $\textcolor{red}{\| f \|_2 \equiv \sqrt{ \sum |c_k|^2 }}$ 
\\[4pt]
(Roth, '54) & 
\\
\hline
&
\\[-3pt]
$\textcolor{red}{\|D_N\|_\infty \gtrsim \log N}$  &  $\textcolor{red}{\| f\|_\infty \gtrsim \sum | c_k |}$ \,\, 
\\[4pt]
(Schmidt, '72; Hal\'asz, '81) & (Sidon, '27)
\\[4pt]
Riesz product: $\prod  \bigl( 1+  c f_k \bigr) $ & Riesz product: $\prod \big(1+ \cos (n_k x +\phi_k)\big)$
\\[4pt]
\hline
&
\\[-3pt]
$\textcolor{red}{\| D_N \|_1 \gtrsim \sqrt{\log N}}$  & $\textcolor{red}{\| f\|_1 \gtrsim \|f\|_2 }$ \,\, 
\\[4pt]
(Hal\'asz, '81) & (Sidon, '30)
\\[4pt]
Riesz product: $\prod  \bigg( 1+  i\cdot \frac{c}{\sqrt{\log N}}  f_k \bigg) $ & Riesz product: $ \prod \bigg(  1+ i \cdot \frac{|c_k|}{\|f\|_2 }   \cos (n_k x +\theta_k)\bigg)$
\\[6pt] \hline
\end{tabular}\vskip1mm
\caption{Discrepancy function and lacunary Fourier series} \label{disclac}
\end{table}

The power of the Riesz product approach in discrepancy problems and the Small Ball Conjecture can be intuitively justified. The maximal values of the discrepancy function (as well as of hyperbolic Haar sums) are achieved on a  very sparse, fractal set. Riesz products are known to capture such sets extremely well. In fact,  $\Psi = 2^{n+1} {\mathbf{1}}_E$, where $E$ is the set on which all the functions $f_k$ are positive, i.e. $\Psi$ defines a uniform measure on the set where the $L^\infty$ norm is achieved. In particular,
$ E $  is essentially the low-discrepancy van der Corput set \cite{MR2817765} if 
 all $\varepsilon_R =1$ (in this case, $f_{(k,n-k)}$ are Rademacher functions).

 Historically, Riesz products were designed to work with lacunary Fourier series, see e.g. \cite{MR1963498}, \cite{Riesz}, \cite{sidon1}, \cite{sidon2}), that is, Fourier series with harmonics supported on lacunary sequences $\{n_k\}$ with $n_{k+1}/n_k > \lambda >1$, e.g., $n_k=2^k$. The terms of such series behave like  independent random variables, which resembles our situation, since the functions $f_{(j,n-j)}$ are  actually  independent. The failure of the product rule explains the loss of independence in higher dimensions (see \cite{MR2731041} for this approach towards  the conjecture). The strong similarity of the two-dimensional Small Ball inequality and Sidon's theorem on lacunary Fourier series \cite{sidon1}
 \begin{equation}
 \bigg\| \sum_{  |R|=2^{-n}} \alpha_R h_R\bigg\|_\infty
\gtrsim 2^{-n}  \sum_{R:\, |R|=2^{-n}} |\alpha_R|     \qquad \textup{ vs. } \qquad \bigg\| \sum_{k} c_k \sin  n_k x \bigg\|_\infty \gtrsim \sum_k |c_k|
 \end{equation}
may be  explained heuristically: the condition $|R|=2^{-n}$ effectively leaves only one free parameter,  and  the supports of Haar functions are dyadic -- thus we obtain a one-parameter system with lacunary frequencies. The similarities between discrepancy estimates, lacunary Fourier series, and the corresponding Riesz product techniques are shown in Table \ref{disclac}.


\section{Recent Results}\label{s.rr}

An improvement of the Small Ball inequality in higher dimensions has been obtained by Bilyk, Lacey, and Vagharshakyan \cite{MR2414745,MR2409170}.

\begin{theorem}\label{t.sb3+} For all dimensions $ d\ge 3$, there is an $ \eta = \eta (d)> c / d ^2 $ so that for all 
integers $ n$ there holds 
\begin{equation*}
2 ^{-n} \sum _{\lvert  R \rvert= 2 ^{-n} } \lvert  \alpha_R  \rvert 
{}\lesssim{}  n ^{\frac{d-1}2 - \eta  }
\Bigl\lVert  \sum _{\lvert  R \rvert\ge 2 ^{-n} } \alpha_R h_R \Bigr\rVert_{\infty }
\end{equation*}
\end{theorem}

\noindent We shall briefly explain some ideas and complications that arise in the  higher-dimensional case.

All simple approaches to these questions are blocked by the dramatic failure of the product rule 
 in dimensions $ d\ge 3$.  This failure, as well as potential remedies, 
was first addressed in the breakthrough paper of J\'ozsef Beck \cite{MR1032337}.   Recall that the product rule breaks when some sides of the dyadic rectangles coincide. There is a whole  range of inequalities which partially compensate for the absence  of the product rule and the presence of coincidences. The simplest of these inequalities is the so-called {\it{Beck gain}}. 

\begin{lemma}[Beck gain]\label{l.SimpleCoincie} In dimensions $ d\ge 3$ there holds 
\begin{equation}\label{e.bg}
\Bigl\lVert  \sum_{ \substack{r\neq s \;:\; \lvert  r\rvert= \lvert  s\rvert=n  \\ r_1 = s_1 }} f _{r} \cdot f _s \Bigr\rVert_{p} 
\lesssim p ^{ d-1} n ^{ \frac{2d - 3}{2}} \,, \qquad 1< p < \infty \,. 
\end{equation}
\end{lemma}

The meaning of this bound can be made clear  by simple parameter counting.  The summation conditions $|r|=|s|=n$ and $r_1=s_1$ ``freeze" 3 parameters. Thus the pair of  vectors $r$ and $s$ has $2d-3$ free parameters, and the estimate says that they behave in an orthogonal fashion, nearly as if we had just applied 
the Littlewood-Paley inequality $ 2d-3$ times. The actual proof is more complicated, of course, 
since the variables in the sum are not free as they are in \eqref{e.slp}.  The paper of Beck \cite{MR1032337} contains a weaker version of the 
lemma above in the case of $ d=3$,  $ p=2$.  The $ L ^{p}$ version is far more useful: the case 
$ d=3$ is in \cite{MR2414745}, and an induction on dimension argument \cite{MR2409170} proves the general case.  

To apply the Riesz product  techniques one has to be able to deal  with longer, more complicated patterns of coincidences. This would require inequalities of the type
\begin{equation}\label{beckgain}
\left\| \sum f_{{r_1}} \cdots f_{{r_k}} \right\|_p \lesssim p^{\alpha M} n^{\frac{M}{2}},
\end{equation}
where the summation is extended over all $k$-tuples of $d$-dimensional integer vectors ${r_1}$, ..., ${r_k}$ with a specified configuration of coincidences and $M$ is the number of free parameters imposed by this configuration,  i.e. the free parameters  should still behave orthogonally even for longer coincidences.  If $k=2$,  this is just \eqref{e.bg}; in \cite{MR2409170} a partial result in this direction is obtained for $k>2$. 

While the breakdown of the product rule  is a feature of the method, there are intrinsic issues that demonstrate  that the higher-dimensional inequality is much more delicate and difficult than the case $d=2$. 
There is no simple closed form for the dual function in this situation.
Indeed, assume that    all $ \lvert  \alpha _R\rvert=1 $. One  then wants to 
show that the sum $ \sum_{R \;:\; \lvert  R\rvert= 2 ^{-n}  } \alpha _R h _R (x) \gtrsim n ^{d/2}$ for some values of  $x$.
But every $ x$ is contained  in many more,namely  $ c n ^{d-1} \gg n ^{d/2}$, rectangles of volume $ 2 ^{-n}$.  
That is, one has to identify a collection of  points which capture  only a very slight disbalance between the number of positive and negative summands. There doesn't seem to be  any canonical way to select such a set of points in the higher-dimensional setting, let alone construct  a function similar to the Riesz product \eqref{psi}, which would be close to uniform measure on such a  set, see \cite{MR2731041}.

\subsection{Other Endpoint Estimates}

The  Small Ball Conjecture provides  supporting evidence for Conjecture \ref{j.D} on the behavior of the $ L ^{\infty }$ norm of 
the discrepancy function in dimensions $ d\ge 3$, $\| D_N\| \gtrsim (\log N)^{d/2}$.  On the other hand, the best known  examples of point sets $ \mathcal P_N$ satisfy 
$ \lVert D_N\rVert_{\infty } \lesssim (\log N) ^{d-1}$.   However, the techniques of the orthogonal function method  cannot prove anything better than the Small Ball inequality.

As we have pointed out repeatedly, the set on which $ D_N$ achieves its $ L ^{\infty }$ norm  is a small set.  Exactly how small 
has been quantified in the two-dimensional setting by  Bilyk, Lacey, Parissis, Vagharshakyan \cite{MR2573600}. 

\begin{theorem}\label{t.blpv} In dimension $ d=2$, for any integer $ N$

(a) for any point set $\mathcal P_N$ with $\#\mathcal P_N = N$, and $ 2 < q < \infty $,  we have  
\begin{equation}
\lVert D_N\rVert_{\operatorname {exp}(L ^{q})}  \gtrsim  (\log N) ^{1- 1/q} ;
\end{equation}

(b)  there exists a set $ \mathcal P_N$ 
(a shifted van der Corput set)  such that for $ 2\le q < \infty $, 
\begin{equation*}
\lVert D_N\rVert_{\operatorname {exp}(L ^{q})}  \lesssim  (\log N) ^{1- 1/q} \,. 
\end{equation*}
\end{theorem}

This theorem is an interpolation between Roth's and Schmidt's bounds in dimension two: when $q=2$ (the subgaussian case) the estimates resembles the $L^2$ behavior, $\sqrt{\log N}$, while as $q$ approaches infinity, the bounds become close to the $L^\infty$ estimate, $\log N$.

The crucial index 
$ q=2$ is the exact limit of  Roth's Theorem:  $ \lVert D_N\rVert_{\operatorname {exp}(L ^{2})}\gtrsim \sqrt{\log N} $ 
by Roth's theorem, and there is an example of $ \mathcal P_N$ for which the reverse inequality holds.  
It is very tempting to speculate that the Orlicz space $ \operatorname {exp}(L^2)$  of subgaussian functions  is the sharp space in all dimensions.  

\begin{conjecture}\label{j.exp} For all dimensions $ d$
\begin{equation*}
\inf _{\mathcal P_N}  \lVert D_N\rVert_{\operatorname {exp}(L ^{2})} \lesssim (\log N) ^{ (d-1)/2} \,. 
\end{equation*}
\end{conjecture}

\noindent This would imply that in the extremal case  the set  $ \{ x \;:\;  D_N (x) \ge (\log N) ^{d/2} \}$ would have measure at most $ N ^{-c}$, 
for some positive $ c$.  We are of course very far from verifying such conjectures, though they can be helpful in devising  potential proof strategies for the main goal -- Conjecture \ref{j.D}.

\end{document}